\documentclass[12pt]{amsart}
\usepackage{amssymb}
\usepackage{amscd}

\newcommand{\Z}{{\mathbb Z}}

\newcommand{\F}{{\mathbb F}}
\newcommand{\R}{{\mathbb R}}

\newcommand{\Fq}{{\F_q}}
\newcommand{\Fqbar}{{\overline{\F}_q}}
\newcommand{\ovF}{{\,\overline{\!F}}}
\newcommand{\alphabar}{\overline{\alpha}}
\newcommand{\OO}{{\mathcal O}}
\newcommand{\PP}{{\mathbb P}}
\newcommand{\C}{{\mathcal C}}
\newcommand{\cF}{{\mathcal F}}
\newcommand{\Q}{{\mathbb Q}}
\newcommand{\isom}{\simeq}
\newcommand{\frakp}{\mathfrak{p}}
\newcommand{\frako}{\mathfrak{o}}
\newcommand{\frakO}{\mathfrak{O}}
\newcommand{\eps}{\varepsilon}

\newcommand{\Drinfeld}{Drinfel'd}
\newcommand{\Vladut}{Vl\u adu\c t}

\def\ly{\fontencoding{U}\fontfamily{lasy}\selectfont}
\newcommand{\guillemetouvert}{{\ly(\kern-0.20em(}}
\newcommand{\guillemetferme}{{\ly)\kern-0.20em)}}

\DeclareMathOperator{\Div}{Div}
\DeclareMathOperator{\End}{End}
\DeclareMathOperator{\Gal}{Gal}

\DeclareMathOperator{\Jac}{Jac}
\DeclareMathOperator{\ord}{ord}

\DeclareMathOperator{\PSL}{PSL}

\theoremstyle{plain}
\newtheorem  {thm}        {Theorem}     [section]
\newtheorem  {lemma}[thm] {Lemma}
\newtheorem  {prop} [thm] {Proposition}
\newtheorem  {cor}  [thm] {Corollary}
\newtheorem* {thm*}       {Theorem}
\newtheorem* {lemma*}     {Lemma}
\newtheorem* {prop*}      {Proposition}

\theoremstyle{definition}

\theoremstyle{remark}
\newtheorem* {remark}     {Remark}

\begin{document}


\title[Curves of every genus]{
  Curves of every genus with many points, II:
  Asymptotically good families
}

\author[Elkies]{Noam D. Elkies}  
\address{
  Department of Mathematics,
  Harvard University,
  Cambridge, MA 02138--2901
}
\email{elkies@math.harvard.edu}

\author[Howe]{\hbox{Everett W. Howe}}    
\address{
  Center for Communications Research,
  4320 Westerra Court,
  San Diego, CA 92121--1967
}
\email{however@alumni.caltech.edu}

\author[Kresch]{\hbox{Andrew Kresch}}
\address{
  Department of Mathematics,
  University of Pennsylvania,
  Philadelphia, PA 19104--6395
}
\email{kresch@math.upenn.edu}

\author[Poonen]{\hbox{Bjorn Poonen}}
\address{
  Department of Mathematics,
  University of California,
  Berkeley, CA 94720--3840
}
\email{poonen@math.berkeley.edu}

\author[Wetherell]{\hbox{Joseph L. Wetherell}}
\address{
  Department of Mathematics,
  University of Southern California,
  Los Angeles, CA 90089--1113
}
\curraddr{
  Center for Communications Research,
  4320 Westerra Court,
  San Diego, CA 92121--1967
}
\email{jlwether@alum.mit.edu}

\author[Zieve]{\hbox{Michael E. Zieve}}
\address{
  Center for Communications Research,
  805 Bunn Drive,
  Princeton, NJ 08540
}
\email{zieve@idaccr.org}

\thanks{The authors thank Jordan Ellenberg, Yasutaka Ihara, and
Felipe Voloch for helpful conversations.
N.\thinspace{}E.\ and B.\thinspace{}P.\ were supported in part by
Packard Fellowships;
in addition, B.\thinspace{}P.\ was supported by NSF Grant DMS-9801104 and
a Sloan Fellowship.
A.\thinspace{}K., J.\thinspace{}W., and M.\thinspace{}Z.\ were
supported in part by NSF Mathematical Sciences Postdoctoral Research Fellowships.
}

\begin{abstract}
We resolve a 1983 question of Serre by constructing
curves with many points of every genus over every finite field.
More precisely, we show that for every prime power $q$ there is a
positive constant $c_q$ with the following property:
for every integer $g\ge 0$, there is a genus-$g$ curve over $\Fq$
with at least $c_q g$ rational points over $\Fq$.
Moreover, we show that there exists a positive constant $d$
such that for every $q$ we can choose $c_q=d\log q$.
We show also that there is a constant $c>0$ such that for every $q$
and every $n>0$, and for every sufficiently large $g$,
there is a genus-$g$ curve over $\Fq$
that has at least $cg/n$ rational points
and whose Jacobian contains
a subgroup of rational points isomorphic to $(\Z/n\Z)^r$ for some $r>cg/n$.
\end{abstract}

\date{23 July 2002}

\maketitle


\section{Introduction}

Let $q$ be a power of a prime $p$ and let $\Fq$ be the field
with $q$ elements.
We will study the function
\begin{equation*}
         N_q(g):=\max \{\#C(\Fq): \text{$C$ is a curve of genus $g$ over $\Fq$} \},
\end{equation*}
where by a {\em curve\/} over a field $k$ we mean a
smooth, projective, geometrically irreducible $1$-dimensional variety over $k$.
In particular, we will be interested in upper and lower bounds
on the asymptotic quantities
\begin{equation*}
    A^{+}(q):=\limsup_{g\rightarrow\infty}N_{q}(g)/g
    \text{\qquad and\qquad}
    A^{-}(q):=\liminf_{g\rightarrow\infty}N_{q}(g)/g.
\end{equation*}
The quantity that we call $A^{+}(q)$ is
denoted $A(q)$ in most current literature.

We begin by reviewing the known bounds on~$A^{+}(q)$.
The upper bound $A^{+}(q)\le 2\sqrt{q}$ follows from
the well-known inequality
\begin{equation*}
        N_q(g) \leq (2\sqrt{q})g + (q+1),
\end{equation*}
which was proved by Weil~\cite{We} in the 1940s.
Serre writes~\cite{Se2} that for a long time Weil's inequality
was considered ``essentiellement \guillemetouvert optimale\guillemetferme.''
The situation changed in 1981, when Manin~\cite{Ma} (following Goppa)
and Ihara~\cite{Ih-81} found bounds on $N_q(g)$ that are significantly
better than the Weil bound when $g$ is large compared to~$q$.
Ihara's idea was refined by \Drinfeld\ and \Vladut~\cite{DV}
to prove that, for fixed~$q$, we have
\begin{equation}
\label{eqn-DV}
        N_q(g)\leq (\sqrt{q}-1+o(1)) g \qquad\text{as $g \rightarrow \infty$,}
\end{equation}
which shows that $A^{+}(q)\le\sqrt{q}-1$.

The problem of finding lower bounds for $A^{+}(q)$ leads to
the problem of producing curves over $\Fq$ with many points.
Serre~\cite{Se,Se-Harvard} used class field towers to
construct infinite sequences of such curves; in this way, he
showed that $A^{+}(q)>0$ for all~$q$.
In fact he proved that $A^{+}(q)>c \log q$ for some $c>0$ not depending on~$q$,
and his proof gave better bounds for nonprime~$q$.
Ihara~\cite{Ih-66, Ih-69, Ih-79, Ih-81} used supersingular points on
Shimura curves to show that if $q$ is a square
we have $A^{+}(q)\geq\sqrt{q}-1$; together
with the \Drinfeld-\Vladut\ bound~(\ref{eqn-DV}), this proves that
$A^{+}(q)=\sqrt{q}-1$ for every square~$q$.
Several other authors proved further results along similar lines.
(For a survey, see the introduction to~\cite{KWZ}.)
But qualitatively the situation today is similar to that in~1983:
\begin{itemize}
  \item There is no nonsquare $q$ for which one knows whether
        or not $A^{+}(q)=\sqrt{q}-1$.
  \item For square $q$, every known sequence $\{C_i\}$ of curves over $\Fq$
        whose genera $g_i$ tend to infinity and for which
        $\#C_i(\Fq)/g_i \rightarrow\sqrt{q}-1$
        has the property that $C_i$ is modular (of elliptic, Shimura, or
        \Drinfeld\ types) for all sufficiently
        large~$i$ (see~\cite{E, E-Drinfeld}).
\end{itemize}

In this paper we consider the problem of finding lower bounds for~$A^{-}(q)$.
To show that $A^{-}(q)$ is greater than some
constant~$c$, one must show that for {\em every\/} sufficiently large $g$
there is a curve of genus $g$ over $\Fq$ with more than $cg$ points.
Unfortunately, the class field tower and modular curve constructions
mentioned above miss infinitely many genera~\cite{CWZ},
so they do not lead immediately to lower bounds for~$A^{-}(q)$.
In fact, before now there was no $q$ for which it was known whether $A^{-}(q)>0$;
furthermore, even the much weaker assertion $\lim_{g\to\infty} N_q(g)=\infty$
was proved only recently~\cite{KWZ}.

The main result of this paper gives an affirmative
answer to a 1983 question of Serre~\cite{Se2}:
\begin{thm}
\label{thm-main}
For every~$q$, we have $A^{-}(q) >  0$.
\end{thm}
We also discuss quantitative refinements of this result. For instance,
at the end of Section~\ref{sec-cft} we show that there is a constant $d>0$
such that  $N_q(g) \ge (d \log q) g$.  In particular, we find that
$A^{-}(q) \ge d\log q$ for all~$q$.
This can be strengthened for square~$q$:
\begin{thm}
\label{thm-mainsquarelog}
We have
$$A^{-}(q)\geq
\begin{cases}
\displaystyle
\frac{\sqrt{q}-1}{2 + (\log 2 / \log q)} &
   \text{if $q$ is an even square{\rm;}} \\
 & \\ 
\displaystyle
\frac{\sqrt{q}-1}{2 + (\log 4 / \log q)} &
   \text{if $q$ is an odd square.}
\end{cases}$$
\end{thm}
\noindent For odd square $q < 417$, Corollary~\ref{cor-smallq}
in Section~\ref{sec-limitations} gives an improved lower bound.

Theorem~\ref{thm-main} shows that for every sufficiently large integer $g$
there exist curves of genus $g$ with many points.
Theorem~\ref{thm-nrank} below shows that we can simultaneously force the
rational points of the Jacobian to contain certain groups.
Suppose $C$ is a curve over a field $k$ and $n>1$ is an integer.
We define the {\em $n$-rank\/} $r_n(C)$ of $C$
to be the largest integer $r$ such that
the group $(\Z/n\Z)^r$ can be embedded in the group of $k$-rational
points of the Jacobian of $C$.
\begin{thm}
\label{thm-nrank}
There exists a constant $c>0$ such that for all integers $n > 1$ and prime
powers $q$ the following statement holds{\/\rm:}
For all sufficiently large integers~$g$, there is a genus-$g$
curve $C$ over\/ $\Fq$ with $r_n(C) > (c/n) g$ and $\#C(\Fq)> (c/n) g$.
\end{thm}

Our proofs of these results involve covers of curves.
A {\em cover\/} of a curve $C$ over a field $k$
is another curve $B$ over $k$ together with a nonconstant
separable map $B\to C$ over~$k$.
Our strategy for proving Theorems~\ref{thm-main}
and~\ref{thm-mainsquarelog}
is to begin with some sequence of curves with many points that
achieve ``enough'' genera,
and then to fill in the missing genera via degree-2 covers.
More specifically, we start with a sequence of curves over $\Fq$ that
have many points and whose genera grow at most exponentially.
Then we show that for every curve $C$ in this sequence,
and for every integer $h$ greater than some constant multiple of the
genus of $C$, there exists a degree-2 cover $B\to C$ over $\Fq$
such that $B$ has genus~$h$.
Either $B$ or its quadratic twist $B'$
will have at least as many $\Fq$-rational points as $C$,
because $\#B(\Fq)+\#B'(\Fq)=2\#C(\Fq)$.
To produce the sequence of curves we use class field towers and Shimura curves
(see Section~\ref{sec-asymptotic}).
To produce the degree-$2$ covers, we use the following result of independent
interest:
\begin{prop}
\label{prop-intro}
Let $C$ be a curve of genus~$g$ over a finite field $k$.  For every integer~$h\ge 4g$,
there exists a genus-$h$ curve $B$ over $k$
that admits a degree-$2$ covering map $B\to C$ defined over~$k$.
\end{prop}
\noindent

In the next section we prove Proposition~\ref{prop-intro},
as well as a generalization of the proposition to curves
over arbitrary fields.
We prove Theorem~\ref{thm-main} and a weak version of
Theorem~\ref{thm-mainsquarelog}
in Sections~\ref{sec-cft} and \ref{sec-shimura}, respectively.
In Section~\ref{sec-largerank} we prove Theorem~\ref{thm-nrank}
and in Section~\ref{sec-constants} we prove the full version
of Theorem~\ref{thm-mainsquarelog}.
We end in Section~\ref{sec-limitations} by discussing limitations on
the results obtainable by our methods.
Our proof of Theorem~\ref{thm-main}
relies on Serre's proof~\cite{Se-Harvard} that $A^{+}(q)>0$ for all~$q$;
since Serre's proof has not been published,
we include a version of it as an appendix.


\section{Degree-2 covers}
\label{sec-doublecover}
In this section we prove Proposition~\ref{prop-intro} and give an
analogous result for curves over infinite fields.

\subsection{Finite fields}
Suppose $C$ is a curve over $\Fq$.
Let $g(C)$ denote the genus of~$C$,
and let $\Fq(C)$ denote the function field of $C$.
For every positive integer $d$, let $n_d(C)$ denote the number
of places of $C$ of degree $d$.
We begin with a lemma
that gives a lower bound on $n_d(C)$ for large enough $d$.

\begin{lemma}
\label{lemma-placecount}
Let $C$ be a curve over\/ $\Fq$ of genus~$g$.
\begin{enumerate}
\item[\rm (i)] For every $d>0$ we have $n_d(C) > (q^d-(6g+3)q^{d/2})/d$.
\item[\rm (ii)] If $g>1$ and $d > 2g$ we have $n_d(C) > 0$.
\item[\rm (iii)] If $q$ is odd and $d > 2g$ then
  $n_d(C)\ge 2^{2g}$,
 with equality if and only if $q=d=3$, $g=1$, and $\#C(\F_3)=6$.
\item[\rm (iv)] Assume $g \ge 900$. If $j$ is a positive integer
and $d \ge \log_q(j \log j+1) + \sqrt{g}$, then $n_d(C) > j$.
\end{enumerate}
\end{lemma}

\begin{proof}
(i)
Clearly $\#C(\F_{q^d})=\sum_{m\mid d} m \cdot n_m(C)$, so
by M\"obius inversion we find that
\begin{equation}
\label{eqn-placeequality}
 d n_d(C) = \sum_{m\mid d} \mu(d/m) \#C(\F_{q^m}).
\end{equation}
Combining the Weil lower bound on $\#C(\F_{q^d})$ with the
Weil upper bounds on the various~$\#C(\F_{q^m})$, we find that
\[
 d n_d(C) \ge q^d + 1 - 2g q^{d/2}
         - \sum_{\stackrel{\scriptstyle m\mid d}{\scriptstyle 1\le m<d}}
                \left( q^m + 1 + 2g q^{m/2} \right).
\]
Using the fact that
$q + q^2 + \dots + q^{\lfloor d/2\rfloor} < 2q^{d/2}$, we see that
\[
\sum_{\stackrel{\scriptstyle m\mid d}{\scriptstyle 1\le m<d}}
                \left( q^m + 1 + 2g q^{m/2} \right)
 <
\sum_{1\le m\le d/2} \left( (2g+1) q^m + 1\right)
 <
(4g+2) q^{d/2} + d/2.
\]
Thus $d n_d(C) > q^d - (6g+2)q^{d/2} - d/2$, and~(i) follows.

(ii)
Suppose $g\ge 2$ and $d\ge 2g+1$.  If $q>2$ then $q^{d/2}\ge q^{g+1/2} > 6g+3$,
so by~(i) we have $n_d(C)>0$.  Suppose $q=2$.  If $g\ge 5$ or if $d\ge g+6$
then once again we have $2^{d/2} > 6g+3$, and (i) shows that $n_d(C)>0$.
We are left to consider six triples $(q,d,g)$.
For each one, we can show that $n_d(C) > 0$ by combining
equation~(\ref{eqn-placeequality})
with the refined Weil bound
$|\#C(\Fq) - q - 1|\le g\lfloor 2\sqrt{q}\rfloor$
(see \cite[Th\'eor\`eme 1]{Se})
and the trivial bound $\#C(\Fq) \ge 0$.

(iii)
We assume $g\ge 1$, since (iii) is well-known when $g=0$
(it amounts to the existence of irreducible
polynomials over $\F_q$ of arbitrary positive degree).
If $q\ge 5$ then (iii) follows from (i), so assume $q=3$.
In this case (iii) again follows from (i) except when $(g,d)$ is
one of $(1,3), (1,4), (2,5)$.  When $(g,d)$ is either $(1,4)$ or
$(2,5)$ then (iii) follows from~(\ref{eqn-placeequality}) and the
Weil bounds.  When $(g,d) = (1,3)$, the statement can be proved by
considering the possible zeta functions of elliptic curves over $\F_3$.

(iv)
Define $f(x) = 2^{\sqrt{x}/2}/2 - (6x+3)$.  A short computation
shows that $f(900)$ is positive and that $f(x)$ is increasing for
$x\ge900$.  We conclude that $(6g+3) \le 2^{\sqrt{g}/2}/2 \le
q^{d/2}/2$; applying~(i), we find that $n_d(C) \ge q^d/2d$.

Note that $d \ge 30$ and that $q^d/2d$ is an increasing function of
$d$ for $d$ in this range.  When $j=1$, this shows that $q^d/2d \ge
2^{30}/60 > j$.  When $j \ge 2$, we find that
$$ \frac{q^d}{2d} \ge
        \frac{(j \log j + 1) q^{30}}{2\log_q(j \log j + 1) + 60} \ge
        \frac{2^{30} j  \log j}{4\log_2 j + 60} \ge
        \frac{2^{30} j  \log j}{64\log_2 j} > j.
$$
In both cases, $n_d(C) > j$, as desired.
\end{proof}

\begin{proof}
[Proof of Proposition~{\rm\ref{prop-intro}}]
The assertion is clear when $g=0$, so we assume $g\ge 1$.

First suppose $q$ is even.
If $\#C(\Fq)>0$ then let $R$ be a rational point of $C$.
Set $m=2h-4g+1$; by Riemann-Roch,
there exists a function $f$ in $\Fq(C)$ having
polar divisor $mR$.
Let $B\to C$ be the cover of $C$
obtained by adjoining a root of $X^2+X+f$ to~$\Fq(C)$.
This cover has degree two; its different is $(m+1)R'$,
where $R'$ is the point of $B$ lying over $R$.
By Hurwitz, the genus of $B$ is $h$, and we are done.
We are left with the case where $\#C(\Fq)=0$.  Since every
genus-$1$ curve over a finite field has a rational point,
we must have~$g\ge 2$. Then
Lemma~\ref{lemma-placecount}(ii) guarantees the existence of
a place $P$ of degree $d:=h-2g+1$.
By Riemann-Roch, there exists $f\in\Fq(C)$ with polar divisor $P$.
If we form $B$ as before, the different of the cover $B\to C$ is $2P'$,
where $P'$ is the place of $B$ lying above $P$.
Again the genus of $B$ is $h$, and we are done.

Now suppose $q$ is odd.
Let $J(\Fq)$ be the group of degree-zero divisor classes on $C$
(or equivalently, the group of $\Fq$-rational points on the Jacobian of $C$).
Let $d=h-2g+1$, and let $T_d$ denote the set of places on
$C$ of degree~$d$.
Consider the map $T_d\to J(\Fq)/2J(\Fq)$ defined by~$P\mapsto [P-P_0]$, where
$P_0$ is any fixed degree-$d$ divisor on~$C$.
Suppose that $\#T_d>\#(J(\Fq)/2J(\Fq))$.  Then there are distinct places
$P,P'\in T_d$ with the same image in~$J(\Fq)/2J(\Fq)$,
so there is a function $f$
on $C$ having divisor~$P-P'+2D$ for some divisor $D$.
Adjoining $\sqrt{f}$ to the function
field of $C$ produces a degree-2 cover $B\to C$ over $\Fq$,
with the genus of $B$ equal to $2g+d-1=h$.
Thus we need only show that $\#T_d>\#(J(\Fq)/2J(\Fq))$.
But $\#(J(\Fq)/2J(\Fq))$ is equal to the number of $2$-torsion points in
$J(\Fq)$, and this is at most $2^{2g}$.
Lemma~\ref{lemma-placecount}(iii) supplies the needed inequality so long
as we are not in the exceptional case where $q=d=3$, $g=1$, and $\#C(\Fq)=6$.
But in that case $\#J(\Fq)=6$, so $\#(J(\Fq)/2J(\Fq))=2$ and the result
follows.
\end{proof}

\begin{remark}
The argument given above for even $q$ can be modified to
show that every curve over a
finite field of characteristic $3$ has Artin-Schreier covers of every
sufficiently large genus.  However, the analogous statement in
characteristic $p>3$ is not true; the degree of the different of an
Artin-Schreier cover $B\to C$ is divisible by~$p-1$, so for any such
cover we have $g(B)\equiv g(C)$ mod~$(p-1)/2$.
\end{remark}

\subsection{Arbitrary fields}
Although we will not need it for our main results,
we now discuss a generalization of
Proposition~\ref{prop-intro} to arbitrary fields.
The {\em index\/} of a curve $C$ is the positive integer $I$
such that the image of the degree map $\Div(C) \rightarrow \Z$ equals~$I\Z$;
equivalently, $I$ is the greatest common divisor
of the degrees of all places of~$C$.
A classical result of F.\thinspace{}K.~Schmidt asserts that curves over finite fields
have index 1 (this follows from Lemma~\ref{lemma-placecount}(ii); see
\cite[Cor.~V.1.11]{St} for a more elementary proof).
Over some infinite fields, higher indices are possible;
for instance, a conic over $\Q$ with no rational points has index~$2$.

\begin{prop}
\label{prop-index}
Let $C$ be a curve of genus $g$ over a field $k$
and let $I$ be its index.
If $B \rightarrow C$ is a degree-$2$ cover of $C$ of genus $h$,
then $h \equiv 1 \pmod{I}$.
Conversely, if $h$ is an integer satisfying $h\ge 4g$
and $h \equiv 1 \pmod{I}$,
then there exists a degree-$2$ cover $B \rightarrow C$
of genus $h$.
\end{prop}

\begin{proof}
The index divides $2g-2$, because $2g-2$ is the degree
of the canonical divisor.
If the characteristic of $k$ is not $2$,
then $k(B) \isom k(C)(\sqrt{f})$ for some $f \in k(C)^*$,
and the total degree of the set of points where $f$ has odd valuation
must be a multiple of $2I$, since the divisor of $f$ has degree zero;
now the Hurwitz formula implies $h \equiv 1 \pmod{I}$.
If $k$ has characteristic $2$,
we have $k(B) \isom k(C)(z)$,
where $z$ satisfies $z^2+z=f$ for some $f \in k(C)$.
A local calculation shows that each point of $B$ occurs
with even multiplicity in the different, so the
degree of the different is divisible by $2I$,
and Hurwitz again implies $h \equiv 1 \pmod{I}$.

Now we prove the converse assertion.
Proposition~\ref{prop-intro} handles the case where $k$ is finite,
so we assume $k$ is infinite.
Let $d=h-2g+1$.
Then $d \equiv 2-2g \equiv 0 \pmod{I}$, so we can choose a divisor $D$ on $C$
of degree $d$.
Since $d \ge 2g+1$, the divisor $D$ is very ample, and therefore determines
an embedding $C \rightarrow \PP^{d-g}$ whose image has degree $d$.
The composition of this embedding with a generically chosen linear projection
to $\PP^1$ is a degree-$d$ map $f\colon C \rightarrow \PP^1$.
The linear projection may be chosen so as not to kill a chosen tangent
vector at some point of $C$, so we may assume the map $f$ induces a separable
extension of function fields.
By composing $f$ with an automorphism of $\PP^1$,
we may further assume that $0$ and $\infty$
are not branch points of $f$.
Now let $B$ be the cover obtained by adjoining to $k(C)$ a root $y$
of $y^2+y=f$ or $y^2=f$, according as the characteristic of $k$
equals $2$ or not.
The different of $B \rightarrow C$ has degree $2d$ in each case,
so by Hurwitz we have $g(B)=2g-1+d=h$.
\end{proof}


\section{Asymptotic lower bounds}
\label{sec-asymptotic}

In this section we prove Theorem~\ref{thm-main} and
a weakened version of Theorem~\ref{thm-mainsquarelog}.
(The full version of Theorem~\ref{thm-mainsquarelog}
will be proved in Section~\ref{sec-mainsquareproof}.)
As we noted in the introduction, Serre used class field towers to exhibit
a positive lower bound on $A^{+}(q)$ for every~$q$,
and Ihara used Shimura curves to produce a better lower bound
for square~$q$.  We will briefly recall these constructions, and include
some necessary complements.  Our theorems will follow from
these constructions and Proposition~\ref{prop-intro}.


\subsection{Preliminaries}
\label{sec-prelims}

We begin by relating Proposition~\ref{prop-intro} to rational points.

\begin{prop}
\label{prop-twists}
Let $C$ be a curve over\/ $\Fq$ of genus~$g$.  For every integer~$h\ge 4g$, we
have $N_q(h)\geq\#C(\Fq)$.
\end{prop}

\begin{proof}
By Proposition~\ref{prop-intro}, there is a genus-$h$ curve $B$
over $\Fq$ that admits a degree-2 covering map $B\to C$ over~$\Fq$.
The result follows, since either $B$ or its quadratic twist over $C$
has at least as many $\Fq$-rational points as does~$C$.
\end{proof}

\begin{cor}
\label{cor-N_q(g)}
If $h \ge 4g$ then $N_q(h)\ge N_q(g)$.\qed
\end{cor}

Now we give a general technique for producing lower
bounds for~$A^{-}(q)$.
We will call a sequence $\C=\{C_0,C_1,\dots\}$ of curves over $\Fq$
{\em ascensive\/} if $\lim_{i\to\infty}g(C_i)=\infty$.
If $\C$ is an ascensive sequence of curves, define
    $$\gamma(\C):=\liminf_{i\to\infty} \frac{\#C_i(\Fq)}{g(C_{i+1})}.$$
Note that we divide
the number of points of each curve in the sequence by the genus of
the {\em next\/} curve.
It is not hard to show that $0\le\gamma(\C)\le A^{+}(q)$.

\begin{prop}
\label{prop-liminf}
If $\C=\{C_0,C_1,\dots\}$ is an ascensive sequence of curves over\/ $\Fq$,
then $A^{-}(q)\geq  \gamma(\C)/4.$
\end{prop}
\begin{proof}
For any $h\ge 4g(C_0)$, let $i$ be the largest integer for which
$h\ge 4g(C_i)$.  Proposition~\ref{prop-twists} implies that~$N_q(h)\geq \#C_i(\Fq)$;
thus,
\begin{equation*}
\frac{N_q(h)}{h} \geq \frac{\#C_i(\Fq)}{4g(C_{i+1})},
\end{equation*}
and the desired result follows.
\end{proof}


\subsection{Class field towers}
\label{sec-cft}

In this section
we prove Theorem~\ref{thm-main} using Proposition~\ref{prop-liminf}.
To do this, we must produce sequences $\C$ with~$\gamma(\C)>0$.
We accomplish this using class field towers.

Let $F$ be the function field of a curve over $\Fq$.
Let $S$ be a nonempty set of rational
places of~$F$, and let $\ell$ be a prime number.
Let $\ovF$ be an algebraic closure of~$F$.
Then the {\em $(S,\ell)$-class field tower of $F$\/} is the sequence
$\cF=\{F_0,F_1,\dots\}$ of function fields over $\Fq$ defined
inductively as follows.
We start with~$F_0=F$.  Then for every integer $i\ge 0$
we let $F_{i+1}$ be the maximal unramified abelian extension of $F_i$
in~$\ovF$,
of degree a power of~$\ell$, in which every place in $S$ splits completely
into rational places.
The tower is called {\em infinite\/} if each extension $F_{i+1}/F_i$
is nontrivial.

The following result demonstrates the utility of curves whose
function fields have infinite $(S,\ell)$-class field towers.

\begin{lemma}
\label{lemma-cft}
Let $C$ be a curve over\/ $\Fq$ of genus~$g(C)>1$, let $S$ be a nonempty
set of rational places of\/ $\Fq(C)$, and let $\ell$ be a prime.
Suppose the $(S,\ell)$-class field tower of\/ $\Fq(C)$ is infinite.
Then there exists an ascensive sequence $\C$ of curves over\/ $\Fq$ such that
$$\gamma(\C)\geq \frac{\#S}{g(C)-1}\cdot\frac{1}{\ell}.$$
\end{lemma}

\begin{proof}
Let $\cF=\{F_0,F_1,\dots\}$ be the $(S,\ell)$-class field
tower of~$\Fq(C)$.  Since each extension $F_{i+1}/F_i$ is Galois of
$\ell$-power degree, there is a refinement of the tower of fields
$$F_0\subset F_1\subset F_2\subset\cdots$$
to a tower of degree-$\ell$ field extensions
$$F_0=\widehat{F}_0\subset \widehat{F}_1\subset\widehat{F}_2\subset\cdots$$
such that each $F_i$ occurs as some~$\widehat{F}_j$.
Clearly $\widehat{F}_n$ is an unramified extension of $F_0$ in which every
place in $S$ splits completely.  Thus, the genus of $\widehat{F}_n$ is
$1+\ell^n (g(C)-1)$, and the number of rational places of $\widehat{F}_n$ is
at least~$\ell^n\cdot\#S$.  Let $C_n$ be a curve over $\Fq$ such that
$\Fq(C_n)\cong\widehat{F}_n$ and let $\C = \{C_0, C_1,\ldots\}$.
Then we have
\begin{equation*}
\frac{\#C_n(\Fq)}{g(C_{n+1})} \geq \frac{\ell^n \cdot \#S}{1+\ell^{n+1} (g(C)-1)},
\end{equation*}
so that $\C$ is an ascensive sequence satisfying the conclusion of the lemma.
\end{proof}

Theorem \ref{thm-main} follows from Lemma \ref{lemma-cft} and
Proposition \ref{prop-liminf}.
Indeed, for
every fixed~$q$, Serre has constructed a curve $C$ over $\Fq$ of genus~$g>1$
and a nonempty set $S$ of rational places of~$\Fq(C)$,
such that the $(S,2)$-class field tower of $C$ is infinite.
(The appendix contains a version of Serre's construction.)

\begin{remark}
In fact, Serre's construction yields $C$ and $S$ with $\#S>c_1 (\log q)^2$
and $g(C)<c_2\log q$, for some positive constants $c_1$ and~$c_2$.
Let $d = \min(c_1/(8c_2), 1/(2c_2))$.
By considering hyperelliptic curves
when $g \le 4c_2 \log q$, and degree-$2$ covers of the curves from Serre's
tower when $g > 4c_2 \log q$, we conclude that $N_q(g) > (d \log q) g$
for all $q$ and~$g$.
This justifies the claim we made in the introduction,
just before the statement of Theorem~\ref{thm-mainsquarelog}.

%
%

For some nonprime $q$, Serre's proof yields better bounds (as does the method of
Zink~\cite{Zi} based on degenerate Shimura surfaces).
If $q$ is a square, one gets even better results by using Shimura curves
(among other methods --- see~\cite{GS1,GS2,GST}),
as we explain in the following section.
\end{remark}


\subsection{Classical modular curves and Shimura curves}
\label{sec-shimura}

In this section we will prove the following weakened version of
Theorem~\ref{thm-mainsquarelog}:

\begin{thm}
\label{thm-mainsquare}
For every square~$q$, we have $A^{-}(q)\geq (\sqrt{q}-1)/4$.
\end{thm}

\begin{proof}
We first sketch a proof in the special case $q=p^2$ using
the classical modular curves~$X_{0}(\ell)$.
If~$q=p^2$, let $\C=\{X_0(\ell_1),X_0(\ell_2),X_0(\ell_3),\dots\}$,
where $\ell_1<\ell_2<\dots$ are the primes distinct from~$p$.
Then, writing $g_i=g(X_0(\ell_i))$ and $N_i=\#X_0(\ell_i)(\Fq)$,
we have
\begin{equation*} N_i\geq (p-1)(\ell_i+1)/12 \text{\qquad and\qquad}
 g_i =\begin{cases}
 (\ell_i-13)/12 & \text{if $\ell_i\equiv 1$~(mod~12)} \\
 \lfloor (\ell_i+1)/12\rfloor & \text{otherwise,}
 \end{cases}
 \end{equation*}
where the formula for $g_i$ is classical, and the bound on $N_i$ holds
because all supersingular points on $X_0(\ell_i)$ are defined over~$\F_{p^2}$
and their number is at least $(p-1)(\ell_i+1)/12$.
%
%
%
We compute
\begin{equation*}
\gamma(\C) = \liminf_{i\to\infty} \frac{N_i}{g_{i+1}}
\ge \liminf_{i\to\infty} \frac{(p-1)(\ell_i+1)}{\ell_{i+1}+1}
= p-1.
\end{equation*}
(In fact, using the \Drinfeld-\Vladut\ bound one can show that
$\gamma(\C) = p - 1$.)
Now Proposition~\ref{prop-liminf} implies that $A^{-}(p^{2})\geq (p-1)/4$.
This proves Theorem~\ref{thm-mainsquare} for the case~$q=p^2$.

A different proof for the case $q=p^2$ can be given using the Igusa-Ihara field
of modular functions of level $n$ over~$\F_{p^2}$, where $n$ ranges over
the positive integers coprime to~$p$.
(This field is the function field of a twist of $X(n)$
defined over $\F_{p^2}$
such that all its supersingular
points are defined over $\F_{p^2}$: an open subset of this twist
parameterizes elliptic curves $E$ with a pairing-respecting
isomorphism from $E[n]$ to the
Galois module $(\Z/n\Z)^2$ on which the $p^2$-Frobenius automorphism
acts as multiplication by $-p$.)
The relevant facts about these fields were announced in \cite{Ih-66} and
proved in \cite[pp.~166--170]{Ih-69} (see also \cite{Ih-75}).

For more general square fields, we use Shimura curves.
Suppose $q = p^{2m}$ where $p$ is prime and $m\ge1$.
Let $F$ be a totally real number field that has a
prime $\frakp$ whose norm is $p^m$.
Let $F_{\frakp}$ be the completion of $F$ at $\frakp$ and
let $\frako$ denote the ring of integers of $F$.
Let $B$ be a quaternion algebra over $F$ whose ramification locus
contains all but one of the infinite places of $F$ and does not
contain $\frakp$;
standard results of class field theory show that there are
infinitely many such~$B$.
As is explained in Ihara's survey article~\cite{Ih-99},
associated to the data $(F, \frakp, B)$ there is a family of
Shimura curves with many points over $\Fq$ ---
see \cite{Sh-70}, \cite{Ih-79}, and \cite{Mo} for further details.
We will show that this family contains enough curves for us to find
a subfamily $\C$, indexed by the positive integers,
for which $\gamma(\C) = \sqrt{q} - 1$.

Let $N_{B/F}$ denote the reduced norm map of the quaternion algebra $B$
and let $\frakO$ be an arbitrary $\frako$-order in $B$.
Let $\frako[\frakp^{-1}]$ denote the ring
obtained by adjoining
to $\frako$ all of the elements of the fractional ideal $\frakp^{-1}$,
and let $\frakO[\frakp^{-1}]=\frakO\otimes_{\frako} \frako[\frakp^{-1}]$.
We define the quaternionic arithmetic group
$$\Gamma:=\{\,\alpha\in \frakO[\frakp^{-1}] : N_{B/F}(\alpha)=1\}/\{\pm 1\},$$
which can be embedded in
$\PSL_2(\R)\times \PSL_2(F_{\frakp})$ as a discrete subgroup.
The group $\Gamma$ contains a torsion-free subgroup $\Gamma_0$
of finite index.
Associated to $\Gamma_0$ there is a curve $X_0$ over $\F_q$
and a set $S_0$ of $\F_q$-places of $X_0$ such that
$\#S_0\ge (\sqrt{q}-1)(g_0-1)$, where $g_0$ is the genus of $X_0$.
To every finite-index congruence subgroup $\Gamma_1\subset\Gamma_0$
there is an associated unramified cover
$X_1\to X_0$ of degree $(\Gamma_0:\Gamma_1)$,
defined over~$\F_q$,
such that all the places of $S_0$ split completely.
Consequently, to obtain a family of $\F_q$-curves ${\mathcal C}$ with
$\gamma({\mathcal C})=\sqrt{q}-1$ it suffices to show that
$\Gamma_0$ has a sequence of finite-index congruence subgroups
such that the index $d_i$
of the $i{\rm th}$ subgroup tends to infinity
and such that $\lim_{i\to \infty} d_{i+1}/d_i=1$.
We will construct such a sequence.

Let $p_1 < p_2 < \cdots$ be the sequence of primes
that lie below degree-$1$ primes $\frakp_i$ of $F$.
Chebotarev's density theorem shows that the sequence $\{p_i\}$ is infinite
and has positive density in the sequence of prime numbers;
in particular, the sequence $(p_{i+1}+1)/(p_i+1)$ approaches $1$ as $i$
approaches infinity.
Throwing away finitely many primes, we may suppose, for each $i$,
that $p_i\ne p$ and that
$\frakO\otimes_{\frako}(\frako/\frakp_i)$ is a central simple algebra.
By Wedderburn's theorem, we may choose an isomorphism
$\frakO[\frakp^{-1}]\otimes_{\frako}(\frako/\frakp_i) \cong M_2(\F_{p_i})$.
Let $U(p_i)\subset \PSL_2(\F_{p_i})$ denote the subgroup of upper-triangular
matrices and define $\Gamma_0(\frakp_i)$ to be the intersection
of $\Gamma_0$ with the subgroup
$$\{\,\alpha\in\frakO[\frakp^{-1}] : N_{B/F}(\alpha)=1\text{ and }
\alphabar\in U(p_i)\,\}/\{\pm 1\}$$
of $\Gamma$, where $\alphabar$ denotes the image of $\alpha$ in
$M_2(\F_{p_i})$ under the chosen isomorphism.

We claim that the index $d_i := [\Gamma_0:\Gamma_0(\frakp_i)]$ is
equal to~$p_i+1$ when $i$ is sufficiently large.
To see this, note that
for each $i$ the map $\alpha\mapsto\alphabar$ is a homomorphism
$\tau_i\colon \Gamma\to \PSL_2(\F_{p_i})$.
A strong approximation theorem \cite[(5.12)]{Sh-64} applies in this case,
and as a consequence $\tau_i$ is surjective for every $i$.
When $p_i > 3$ the group $\PSL_2(\F_{p_i})$ is simple with order divisible by $p_i$
and hence has no subgroups of index less than~$p_i$.
Thus, for $i$ sufficiently large,
the restriction of $\tau_i$ to $\Gamma_0$ is surjective as well.
It follows that the index $d_i$ is equal to $p_i+1$, as claimed.

The sequence $d_{i+1}/d_i$ approaches $1$.
As we noted above,
this is enough to prove Theorem~\ref{thm-mainsquare} for the field~$\Fq$.
\end{proof}


\section{Curves with many points and large $n$-rank}
\label{sec-largerank}

In our proof of Theorem~\ref{thm-main} we started with an ascensive
sequence of curves with many points and used double covers to produce
curves of every sufficiently large genus with many points.
In this section we will use the same basic method to prove
Theorem~\ref{thm-nrank}: We will start with an ascensive sequence
of curves with many points and large $n$-rank, and use double covers
to produce curves with the same properties in every sufficiently large genus.

Let us begin by introducing some notation.  Suppose
$\C = \{C_0, C_1, \ldots\}$ is an ascensive sequence
of curves over~$\Fq$.  We define $\rho_n(\C)$ (for every
integer $n>1$) and $\beta(\C)$ by
$$\rho_n(\C) := \liminf_{i\to\infty}\frac{r_n(C_i)}{g(C_{i+1})}
\text{\qquad and\qquad}
\beta(\C) := \liminf_{i\to\infty} \frac{g(C_i)}{g(C_{i+1})}.$$
Note that in each expression we are dividing an
invariant of one curve in the sequence by the genus of the next.

We will prove two lemmas.  The first shows how we can obtain
ascensive sequences having positive values of~$\rho_n$.
\begin{lemma}
\label{lemma-r2rb}
Suppose $\C$ is an ascensive sequence over $\Fq$.
Then for every integer $n>1$ there is an
ascensive sequence $\C_n$ over $\Fq$ such that
$$\gamma(\C_n) \ge \frac{\min(\gamma(\C),\beta(\C))}{ 7 n }
\text{\qquad and\qquad}
\rho_n(\C_n)   \ge \frac{\min(\gamma(\C),\beta(\C))}{ 7 n }.$$
\end{lemma}

The second lemma shows how we can get curves with many points
and large $n$-rank from the sequences produced by Lemma~\ref{lemma-r2rb}.
\begin{lemma}
\label{lemma-rb2all}
Suppose $\C$ is an ascensive sequence over $\Fq$.
Then for every $\eps>0$ the following statement holds{\/\rm:}
For every sufficiently large $h$ there is a curve $B$
of genus $h$ such that
$$\#B(\Fq) > (1-\eps) \gamma(\C) h / 4
\text{\qquad and\qquad}
r_n(B)   > (1-\eps) \rho_n(\C) h / 4.$$
\end{lemma}

Our refined versions of Theorem~\ref{thm-main} show that
there is a constant $c$ such that for every $q$
there are ascensive sequences $\C$ over $\Fq$ with
$\beta(\C) = 1$ and $\gamma(\C) > c$.
To prove Theorem~\ref{thm-nrank} we need merely
apply Lemma~\ref{lemma-r2rb} to such a $\C$ and then
apply Lemma~\ref{lemma-rb2all} to the resulting sequences~$\C_n$.

\begin{proof}[Proof of Lemma~{\rm\ref{lemma-r2rb}}]
We may assume that $\C = \{C_0, C_1, \ldots\}$ is an ascensive sequence
over $\Fq$ with $\gamma(\C) > 0$ and $\beta(\C) > 0$.
By eliminating a finite number of terms from $\C$, if necessary,
we may assume that $g(C_i) > 1$ and $\#C_i(\Fq) > 0$ for all~$i$.
Now suppose we are given a positive integer $n$.  Our goal
is to define a sequence $\C_n = \{B_0, B_1, \ldots\}$ satisfying
the conclusions of the lemma. For every $i$, we define $B_i$
as follows.

Let $g = g(C_i)$, let $d = \min(\#C_i(\Fq),g)$, and let $S$
be a subset of $C_i(\Fq)$ with $\#S = d$.
Let $Q$ be a place of $C_i$ of degree $2g + 1$; such a
place exists by Lemma~\ref{lemma-placecount}(ii).
The Riemann-Roch theorem, applied to the divisor $Q$,
shows that there is a function
$f_Q \in \Fq(C_i)$ whose polar divisor is $Q$.
For every $P\in S$, Riemann-Roch applied to $Q$ and $P + Q$ shows
that there is a function $f_P\in \Fq(C_i)$
with simple poles at $P$ and $Q$ and no other poles.  Let
$$f := \sum_{P \in S} f_P,$$
and if $f$ has no pole at $Q$ then replace $f$ with $f + f_Q$.
Then $f$ is a function of degree $d + 2g + 1$
with a simple pole at every $P\in S$.

Write $n = p^e r$ where $p$ is the characteristic of $\Fq$ and $(r,q) = 1$.
Let $B'$ be the curve obtained by
adjoining a root of $y^r = f$ to $\Fq(C_i)$.
If $e=0$, let $B_i=B'$;
otherwise let $B_i$ be the curve obtained by adjoining a root
of $z^{p^e} - z = f$ to $\Fq(B')$.
We note several facts about the curve $B_i$:
First, if we apply the Hurwitz formula to the covers
$B'\to C_i$ and $B_i\to B'$ (the formulas from~\cite[Cor.~III.7.4]{St}
and~\cite[Prop.~III.7.10(e)]{St} are helpful here)
and apply some easy estimates that depend on the
fact that $d\le g$, we find that $g(B_i) < 7 n g$.
%
%
%
%
%
Second, the cover $B_i\to C_i$ is totally ramified at every
point in $S$, so $\#B_i(\Fq) \ge d$.  Third, we have
$r_n(B_i)\ge d - 2$, as the following argument shows.

Let $\pi$ be the natural map from $B'$ to $C_i$.
Let $G$ denote the group of degree-zero divisors of $B'$
supported on $S' := \pi^{-1}(S)$,
and let $G_0$ be the subgroup of $G$ consisting of those degree-zero
divisors $D$ with
$\ord_P D \equiv \ord_{P'} D \pmod r$ for every $P, P'\in S'$.
Note that $G$ is a free abelian group of rank~$d-1$, that
$G_0\supseteq rG$, and that $G_0 / rG$ is naturally isomorphic
to a subgroup of $\Z/r\Z$.

We claim that the kernel of the natural map
$$G \rightarrow (\Jac B')(\Fq)/\pi^\ast((\Jac C_i)(\Fq))$$
contains $rG$ and is contained in $G_0$.
Clearly $rG$ is in the kernel, since if $D \in G$ we have
$rD = \pi^*\pi_*D$.
On the other hand, suppose that $D \in G$ is in the kernel.
Then $D = (h)+\pi^* E$ for some $h \in \Fq(B')$ and divisor $E$ on $C_i$.
The function $h$ is preserved up to scalar multiple by the action
of $\Gal(B'/C_i)$ (by which we mean the Galois group of the covering
of curves over $\Fqbar$),
and writing $h$ in terms of the basis $1,y,\dots,y^{r-1}$
shows that $h \in y^j \Fq(C_i)^*$ for some~$j$.
It follows that the multiplicities of the points of $S'$ in $D$
are all congruent to each other modulo $r$,
so $D$ is contained in $G_0$.
This proves the claim.

It follows from the previous paragraph that the image
of $G$ in the quotient group $$(\Jac B')(\Fq)/\pi^\ast((\Jac C_i)(\Fq))$$
contains $(\Z/r\Z)^{d-2}$.
Group theory then shows
that $(\Jac B')(\Fq)$ also contains $(\Z/r\Z)^{d-2}$.
But since the natural map $B_i\to B'$ has degree coprime to~$r$,
the kernel of the natural map $\Jac B'\to \Jac B_i$ has order coprime
to~$r$, 
so $(\Jac B_i)(\Fq)$ also contains a copy of $(\Z/r\Z)^{d-2}$.

The same proof shows that $(\Jac B_i)(\Fq)$ contains $(\Z/p^e\Z)^{d-2}$.
(The only difference in the argument is that now when
the action of $\Gal(B_i/B')$ on the analogous $h$ is by scalar
multiples, the scalars must be $p^e$-th roots of $1$, hence equal to $1$,
so $h$ comes from~$\Fq(B')$.)
Combining this with the previous paragraph
and applying the Chinese Remainder Theorem,
we find $r_n(B_i) \ge d-2$, as we had noted above.

Combining the several facts we have noted about the curve $B_i$,
we find that for every $i\ge 0$ we have
$$\frac{\#B_i(\Fq)}{g(B_{i+1})} >
      \frac{\min(\#C_i(\Fq),g(C_i))}{7n g(C_{i+1})}
\text{\qquad and\qquad}
\frac{r_n(B_i)}{g(B_{i+1})}     >
      \frac{-2 + \min(\#C_i(\Fq),g(C_i))}{7n g(C_{i+1})}.$$
It follows that
$$\gamma(\C_n)  \ge \frac{\min(\gamma(\C),\beta(\C))}{7n}
\text{\qquad and\qquad}
\rho_n(\C_n)    \ge \frac{\min(\gamma(\C),\beta(\C))}{7n}.$$
\end{proof}

The proof of the second lemma is almost identical to
the proof of Proposition~\ref{prop-liminf}.

\begin{proof}[Proof of Lemma~{\rm\ref{lemma-rb2all}}]
Suppose we are given an ascensive sequence $\C = \{ C_0, C_1, \ldots\}$
and an $\eps>0$ as in the statement of the lemma.
We may certainly assume that $\eps<1$.
Then there is an integer $j$ such that for all $i > j$ we have
$$\#C_i(\Fq)  \ge (1-\eps) \gamma(\C) g(C_{i+1})
\text{\qquad and\qquad}
  r_n(C_i)    \ge (1-\eps) \rho_n(\C) g(C_{i+1}).$$
Now suppose we are given an integer $h \ge 4g(C_{j+1})$.
Let $i$ be the largest integer such that $h \ge 4g(C_i)$.
Proposition~\ref{prop-intro} shows that there is a curve $B$ of
genus~$h$ that is a double cover of~$C_i$.
By replacing $B$ with its quadratic twist over $C_i$,
if necessary, we may assume that $\#B(\Fq) \ge \#C_i(\Fq)$.
Also, the Hurwitz formula shows that the double cover $B\to C_i$
must be totally ramified at some point,
so class field theory for curves shows that
the natural map $(\Jac B)(\Fq)\to (\Jac C_i)(\Fq)$
is surjective.  It follows that $r_n(B) \ge r_n(C_i)$.
Since we also have $h < 4 g(C_{i+1})$, we find that
$$\#B(\Fq) \ge \#C_i(\Fq) \ge (1-\eps) \gamma(\C) g(C_{i+1})
           > (1-\eps) \gamma(\C) h / 4.$$
Similarly, we have
$$r_n(B) \ge r_n(C_i) \ge (1-\eps) \rho_n(\C) g(C_{i+1})
           > (1-\eps) \rho_n(\C) h / 4,$$
as we were to show.
\end{proof}

\begin{remark}
The reader who objects to the intrusion of class field theory 
into the proof of Lemma~\ref{lemma-rb2all} can remove it at the
expense of slightly weakening the lemma.
Consider the double cover $\pi\colon B\to C_i$ introduced at the end of
the proof.  We have $\pi_*\pi^* = 2$ on $\Jac C_i$, so 
$r_n(B) \ge r_{2n}(C_i)$.  Thus, {\em without\/} class field
theory, we can show that for every sufficiently large $h$ there
is a curve $B$ of genus $h$ such that 
$\#B(\Fq) > (1-\eps) \gamma(\C) h / 4$
and
$r_n(B)   > (1-\eps) \rho_{2n}(\C) h / 4$.
\end{remark}

Theorem~\ref{thm-nrank} has the following consequence for class
field towers:

\begin{prop}
\label{prop-classfieldtowerofeverygenus}
Fix a prime power $q$ and a prime $\ell$.
For $g \gg 0$, there exists a genus-$g$ curve $C$ over\/ $\Fq$
with $P \in C(\Fq)$
such that the $(\{P\},\ell)$-class field tower of\/ $\Fq(C)$ is infinite.
\end{prop}

\begin{proof}
Suppose $g$ is at least $5\ell/c$, where $c$ is the constant that
appears in Theorem~\ref{thm-nrank}.  Then Theorem~\ref{thm-nrank}
shows that there is a genus-$g$ curve $C$ over $\Fq$ whose $n$-rank is
greater than~$5$ and that has at least $5$ rational points.
Pick $P \in C(\Fq)$ and embed $C$ in its Jacobian $J$ using $P$ as the
base point.  Let $F \in \End J$ be the $q$-power Frobenius endomorphism
and let $D$ be the inverse image of $C$ under $(1-F)\colon J \rightarrow J$.
Class field theory for curves shows that $D$ is an irreducible curve
and that $D \rightarrow C$ is a unramified Galois covering with Galois 
group $J(\Fq)$, and clearly $P$ splits completely in this covering.
By Galois theory, there exists a subcovering $E \rightarrow C$
with Galois group $(\Z/\ell\Z)^r$, where $r = r_n(C) >5$.
The function field analogue of the Golod-Shafarevich criterion
(the lemma in our appendix, with $S=\{P\}$)
implies that the $(\{P\},\ell)$-class field tower of $\Fq(C)$ is infinite.
\end{proof}


\section{Quantitative refinements}
\label{sec-constants}

Our arguments in Section~\ref{sec-asymptotic} required asymptotic
results about special sequences of curves, but the tool we applied
to these results --- Proposition~\ref{prop-intro} ---
makes no use of any special properties these curves might have.
Indeed, the proposition makes reference only to the genera of the
curves it mentions.
In this section we will show that we can produce
better bounds on $A^{-}(q)$ by replacing Proposition~\ref{prop-intro}
with sharper asymptotic results.  At the end of this section we
prove Theorem~\ref{thm-mainsquarelog}.

\subsection{Bounding data}
\label{sec-boundingdata}

Let $\C=\{C_0, C_1,\ldots\}$ be an ascensive sequence of curves over $\Fq$
and let $H$ and $M$ be real numbers. We will say that the pair $(H,M)$ is
{\em bounding data for $\C$ {\rm(}over $\Fq${\rm)}\/} if
for every $\eps>0$ there is an $L_\eps$ such that
following statement holds:
\begin{itemize}
\item[] If $i>L_\eps$ and $h$ is an
integer with $h > (H+\eps) g(C_i)$, then there exists
a genus-$h$ curve $B$ over $\Fq$ such that $\#B(\Fq) \ge (M-\eps) \#C_i(\Fq)$.
\end{itemize}
Bounding data for a sequence $\C$ can be used to give a lower bound
on $A^-(\C)$, as the following proposition shows.

\begin{prop}
\label{prop-boundingdata}
If $(H,M)$ is bounding data for an ascensive sequence
$\C$ over\/ $\Fq$, then
$A^{-}(q)\geq  \gamma(\C)\cdot M/H.$
\end{prop}
\begin{proof}
This is a generalization of Proposition~\ref{prop-liminf}, and the
proof is essentially identical to the proof of that result.
Indeed, Proposition~\ref{prop-liminf} amounts to the observation that
$(4,1)$ is bounding data for every ascensive $\C$,
which follows from Proposition~\ref{prop-intro}.
\end{proof}

\subsection{Improved bounds on $A^-(q)$ for square $q$}
\label{sec-mainsquareproof}

In this section we will prove Theorem~\ref{thm-mainsquarelog} by showing that
every ascensive sequence over $\Fq$ has bounding data of the form $(H,1)$ for
an appropriate value of $H$.  To show that $(H,1)$ is bounding data for $\C$,
we must show that for every $\eps>0$, for every sufficiently large $i$,
if $h > (H+\eps)g(C_i)$ then there is a genus-$h$ curve $B$ with
$\#B(\Fq)\ge (1-\eps)\#C_i(\Fq)$.
In this section we will restrict ourselves to the case where $B$ is a
degree-$2$ cover of~$C_i$.
This restriction entails that $h\ge 2g(C_i)-1$,
so this line of argument cannot possibly work for values of $H$ less than $2$.

If $\C = \{C_0, C_1,\ldots\}$ is an
ascensive sequence over $\Fq$, define
    $$R_2(\C) := \limsup_{i\to\infty} \frac{r_2(C_i)}{g(C_i)},$$
where, as before, $r_2(C_i)$ denotes the $2$-rank of $(\Jac C_i)(\Fq)$.
(Note that we are dividing the $2$-rank of each curve in the sequence
by {\em its own\/} genus, and {\em not\/} the genus of the next curve!)
Also define
    $$H_\C := 2 + R_2(\C)\frac{\log 2}{\log q}.$$

\begin{prop}
\label{prop-mainsquare}
If $\C$ is an ascensive sequence over $\Fq$ then $(H_\C,1)$ is
bounding data for~$\C$.
\end{prop}

\begin{proof}
Our proof depends on the parity of $q$.
First suppose that $q$ is odd.
Let $C$ be a genus-$g$ curve over $\Fq$.
Let $J$ be the Jacobian of $C$.
Let $j=\#(J(\Fq)/2J(\Fq))$.
Lemma~\ref{lemma-placecount}(iv) shows that
if $g \ge 900$ and $d \ge \log_q(j \log j + 1) + \sqrt{g}$,
then $n_d(C) > j$.
Since $\log j \le \log(2^{2g})=O(g)$,
the condition on $d$ can be rewritten as
    $$d \ge \log_q j + o(g) = r_2(C) \frac{\log 2}{\log q} + o(g)$$
as $g \rightarrow \infty$.
Proposition~\ref{prop-mainsquare} (for odd $q$) follows
by applying the argument in the final paragraph
of the proof of Proposition~\ref{prop-intro}.

Now suppose $q$ is even.  Fix a curve $C/\Fq$ of genus~$g$ and let $d$ be
an arbitrary integer with $d \ge d(C)$, where
    $$d(C) := \max\left( r_2(C)\frac{\log 2}{\log q},
            \frac{2 \log(6g+3)}{\log q}\right).$$
We see from Lemma~\ref{lemma-placecount}(i) that $C$ has a place
of degree~$d$, and it follows from Lemma~\ref{lemma-qeven} (below) that there
exists a degree-$2$ cover $B\to C$ with $g(B) = 2g + d - 1$.  As $C$ ranges
through the curves in $\C$, we have
$$\limsup \frac{d(C)}{g(C)} = R_2(\C)\frac{\log 2}{\log q}.$$
\end{proof}

\begin{lemma}
\label{lemma-qeven}
Let $q$ be a power of $2$, and let $C$ be a curve over $\Fq$.  Suppose $C$
has a place of degree $d$, where $q^d \ge J[2](\Fq)$.
Then there exists a degree-$2$ cover $B\to C$ such that $g(B) = 2g(C) + d - 1.$
\end{lemma}

\begin{proof}
Let $\frakp$ be a place of degree $d$ on $C$ and let $P$ be a geometric point
of $C$ lying in $\frakp$, so that $\frakp$ consists of the $d$ conjugates of $P$.
Let $u$ be a uniformizing parameter at $P$.  For every integer $m\ge 0$
let $S_m$ be the additive group of functions in $\Fq(C)$ that have no poles
outside $\frakp$ and whose polar expansions at $P$ are of the form
$$c_m u^{-2^m} + c_{m-1} u^{-2^{m-1}} + \cdots + c_1 u^{-2} + c_0 u^{-1},$$
where the $c_i$ are elements of $\Fqbar$ and where $c_m$ is not necessarily nonzero.
A standard Riemann-Roch argument shows that $S_m$ contains nonconstant functions
when \mbox{$m \gg0$}, and in fact it is not hard to see that $\#S_m = q^d\#S_{m-1}$
for $m \gg0$.  For $m>0$ let $T_m \subseteq S_m$ be the subgroup
$$T_m := \{g^2 + g + c : g \in S_{m-1}, c\in \Fq\}.$$
Then $\#T_m = \#S_{m-1}$, so for $m\gg0$ we have $\#(S_m/T_m) = q^d$.

Artin-Schreier theory shows the following:
\begin{enumerate}
\item{If $f\in S_m\setminus T_m$ then the degree-$2$ cover of $C$ obtained by adjoining
      a root of $y^2 - y - f$ to $\Fq(C)$ is geometrically irreducible and is
      ramified at most at $\frakp$.}
\item{If $f,g \in S_m\setminus T_m$ then the degree-$2$ covers
      obtained as in~(1) from
      $f$ and $g$ are geometrically isomorphic
      (as curves {\em equipped with their maps to $C$\/})
      if and only if $f - g \in T_m$.}
\end{enumerate}
Thus, from the set $S_m$ we obtain $-1 + \#(S_m/T_m)$ geometrically nonisomorphic
irreducible degree-$2$ covers of $C$.  By taking $m$ large enough, we find that
there are at least $q^d - 1$ such covers.  But the number of geometrically
distinct {\em unramified\/} degree-$2$ covers of $C$ is equal to $\#J[2](\Fq) - 1$,
so there must be a {\em ramified\/} degree-$2$ cover $B\to C$ over $\Fq$ coming
from $S_m$, and it must be ramified exactly at $\frakp$.
The Hurwitz formula then shows that $g(B) = 2g(C) + d - 1$.
\end{proof}

\begin{proof}[Proof of Theorem {\rm\ref{thm-mainsquarelog}}]
We showed in Section~\ref{sec-shimura} that for every square $q$ there exists
an ascensive sequence $\C$ of curves over $\Fq$ with $\gamma(\C) = \sqrt{q} - 1$.
In light of Propositions~\ref{prop-boundingdata} and~\ref{prop-mainsquare},
to prove Theorem~\ref {thm-mainsquarelog} it will be enough for us to show that
$$R_2(\C)\le
\begin{cases}
1 & \text{when $q$ is even;}\\
2 & \text{when $q$ is odd.}
\end{cases}$$
But these inequalities follow from the fact that for every curve $C$ over $\Fq$
we have
$$r_2(C)\le
\begin{cases}
g(C)  & \text{when $q$ is even;}\\
2g(C) & \text{when $q$ is odd.}
\end{cases}$$
\end{proof}


\section{Variants of our argument}
\label{sec-limitations}

Section~\ref{sec-asymptotic} gave a lower bound on $A^-(q)$
in terms of the value of $\gamma(\C)$ of an ascensive sequence.
Section~\ref{sec-boundingdata} gave an improved bound
taking into account also the size of $q$,
but the full power of Proposition~\ref{prop-mainsquare}
was not used, because we used only the trivial upper bound on $R_2(\C)$.
It would be extremely interesting to obtain
better bounds on $R_2(\C)$ when $\C$ is one of the sequences of modular
curves or Shimura curves used in Section~\ref{sec-shimura}.
Ideally, one would like to show that $R_2(\C)=0$ for such an
ascensive sequence, so that $(2,1)$ would be bounding data for $\C$.
Theorem~\ref{thm-nrank} shows, however,
that $R_2(\C)$ can be positive, even for sequences of curves with
many points.  Thus, proving that $R_2(\C)=0$ for a sequence of Shimura curves
would require using special properties of the curves.

It is possible to improve the $4g$ in Proposition~\ref{prop-intro}
if we require only that the proposition apply to curves of sufficiently
large genus.  Indeed, the proof of Proposition~\ref{prop-mainsquare} shows
that we can replace the $4g$ with $(2 + (\log 2 / \log q))g$ if $q$ is even
and with $(2 + (\log 4 / \log q))g$ if $q$ is odd.  But the
following argument limits the possible further
improvements we might hope to obtain:
Given $m$ and $q$, there exists a curve $C$ over $\Fq$
such that every place has degree at least $m$.
(For instance, nonsingular plane curves with this property exist~\cite{Po}.)
It then follows from Proposition~\ref{prop-intro}
that such $C$ exist in every sufficiently large genus.
An argument similar to that proving Proposition~\ref{prop-index}
shows that such $C$ do not admit degree-2 covers of any
genus $h < 2g+2m-1$ except possibly for unramified degree-$2$ covers
of genus $h=2g-1$.
In particular, we cannot improve the $4g$ in Proposition~\ref{prop-intro}
to $2g+s$ for any fixed constant $s$, even if we ask only that
the proposition hold for $g\gg0$.
We do not know whether there always exists a cover of genus $2g+o(g)$.

It is conceivable that $A^-(q) = \sqrt{q} - 1$.
In order to prove this
using our methods, we would need an ascensive sequence $\C$ with
$\gamma(\C) = \sqrt{q} - 1$ for which we could provide bounding
data $(H,M)$ with $H/M = 1$.  But the smallest ratio $H/M$ that we
can possibly attain using degree-$2$ covers $B\to C$ is~$2$, unless we
can make $M>1$ by forcing {\em more\/} than half of the points in $C(\Fq)$
to split.
For $q<207$, it is in fact possible to show
in this way that every ascensive sequence over $\Fq$ has
bounding data $(H,M)$ with $H/M < 2$;
plugging this into Proposition~\ref{prop-boundingdata}
improves the lower bounds on $A^-(q)$.
We give the argument below, beginning with the following variant
of Proposition~\ref{prop-intro}.

\begin{prop}
\label{prop-introplus}
Let $q$ be an odd prime power and let $m = (\log 2 / \log q)(\sqrt{q}+1)$.
Fix a positive real $\eps$.  Then for all $g\gg0$, if $C/\Fq$ is a
curve of genus $g$, and if $h$ is an integer with $h \ge (2 + m + \eps)g$,
then there exists a genus-$h$ curve $B/\Fq$ that admits a
degree-$2$ covering map $B\to C$ in which every rational place of $C$ splits.
\end{prop}

\begin{proof}
Let $d = h - 2g + 1 > (m + \eps) g$.
Thus $d>1$ when $g \gg 0$.
As in the proof of Proposition~\ref{prop-intro},
let $T_d$ denote the set of places on $C$ of degree~$d$.
By Lemma~\ref{lemma-placecount}(i), we
have $\#T_d \ge q^d / (2d)$ when $g\gg0$.
Thus $g \gg 0$ implies $\#T_d \ge q^{(m + \eps/2)g}$.
Consider the map $T_d\to J(\Fq)/2J(\Fq)$ defined by $P\mapsto [P-P_0]$,
where $P_0$ is any fixed degree-$d$ divisor on~$C$.
Since $\#J(\Fq)/2J(\Fq) \le 2^{2g}$,
there exists $x \in J(\Fq)/2J(\Fq)$ having
at least $q^{(m + \eps/2)g} / 2^{2g}$ preimages in $T_d$.
Fix one of these preimages $P_1$.
Then for each preimage $P$ of $x$,
fix $f_P \in \Fq(C)$ having divisor $P - P_1 + 2D$ for some~$D$.

Let $K=\Fq(C)$.
For each $Q\in C(\Fq)$, let $K_Q$ be the completion of $K$ at $Q$,
and let $\OO_Q \subset K_Q$ be the valuation ring.
Since $d>1$, $f_P$ has even valuation at each $Q \in C(\Fq)$.
Hence the $f_P$ map into the group
    $$S := \prod_{Q \in C(\Fq)} \frac{\OO_Q^* K_Q^{*2}}{K_Q^{*2}}$$
of order $2^{\#C(\Fq)}$.
The \Drinfeld-\Vladut\ bound gives
$\#C(\Fq) < (\sqrt{q} - 1 + \eps/2)g$ for $g \gg 0$,
but then the definition of $m$ implies
    $$\frac{q^{(m + \eps/2)g}}{2^{2g}}
    > 2^{(\sqrt{q} - 1 + \eps/2)g} > \#S,$$
so there exist distinct $P,P' \in T_d$ in the preimage of $x$
such that $f_P$ and $f_{P'}$ have the same image in $S$.
Let $e = f_P / f_{P'}$.

Let $B$ be the curve with $\Fq(B) \isom \Fq(C)(\sqrt{e})$.
The divisor of $e$ has the form $P - P' + 2D$ with $P,P' \in T_d$,
so Hurwitz shows that $B$ has genus $h$.
By construction, $e \in K_Q^{*2}$ for all $Q \in C(\Fq)$,
so every rational point of $C$ splits in the degree-$2$ cover $B\to C$.
\end{proof}

\begin{cor}
\label{cor-smallq}
Let $q$ be an odd prime power and let $m = (\log 2 / \log q)(\sqrt{q}+1)$.
Let $\C$ be an ascensive sequence over $\Fq$.
Then $A^-(q) \ge \gamma(\C) / (1 + m/2)$.
\end{cor}

\begin{proof}
Proposition~\ref{prop-introplus} shows
that $\C$ has bounding data $(2+m,2)$.
Now apply Proposition~\ref{prop-boundingdata}.
\end{proof}

We have $1+m/2<2$ for all odd prime powers $q<207$,
and $1+m/2 < 2 + \log 4/\log q$
(the trivial upper bound on $H_\C$)
for all odd prime powers $q<417$.
We showed in Section~\ref{sec-shimura} that, for every square $q$,
there are ascensive sequences $\C$ over $\Fq$ with $\gamma(\C)=\sqrt{q}-1$.
Thus, for odd square $q$ we have $A^-(q)\ge (\sqrt{q}-1)/(1+m/2)$, which
improves the bound in Theorem~\ref{thm-mainsquarelog} if $q<417$.
For example, if $q = 9$ then $1 + m/2 < 1.631$,
and the ascensive sequence $\{X_0(\ell)\}$ over $\F_9$ gives
$A^-(9) > A^+(9) / 1.631 =  (\sqrt{9}-1) / 1.631 > 1.226$.


\appendix
\section*{Appendix: class field towers}
In this appendix we present a version of Serre's proof that $A^{+}(q)>0$
for all~$q$.  Previously Serre's proof has appeared only in the
unpublished lecture notes~\cite{Se-Harvard}.  The version below uses
an idea from~\cite{LM} at one step to simplify the argument.

As in Section~\ref{sec-cft}, it suffices to show that for every $q$
there is a curve $C$ over $\Fq$ of genus $g>1$ such that for some nonempty
set $S$ of rational places on~$\Fq(C)$, the $(S,2)$-class
field tower of $\Fq(C)$ is infinite.  To do this, we use a
function field analogue of the Golod-Shafarevich criterion:
\begin{lemma*}
Let $C$ be a curve over\/ $\Fq$, let $\ell$ be a prime number, and
let $S$ be a nonempty set of rational places of\/ $\Fq(C)$.  Suppose there
is an unramified Galois extension $F/\Fq(C)$,
with Galois group $(\Z/\ell\Z)^r$,
such that every place in $S$ splits completely
into\/ $\Fq$-rational places of~$F$.
If $r\geq 2$ and
\begin{equation*}
\frac{(r-2)^2}{4}\geq \begin{cases}
1+\#S & \text{if $\ell\mid (q-1)$;}\\
\#S   & \text{otherwise,}
\end{cases}
\end{equation*}
then the $(S,\ell)$-class field tower of $\Fq(C)$ is infinite.
\end{lemma*}

\begin{proof}
See~\cite{Sc}.
\end{proof}

\begin{thm*}[Serre]
We have $A^+(q)>0$ for all prime powers $q$.
In fact, there exists $c>0$ such that $A^+(q)> c \log q$ for all $q$.
\end{thm*}

\begin{proof}
First suppose $q$ is odd.  Let $f=f_1f_2\dots f_6$, where $f_1,f_2,\dots,f_6$
are distinct monic irreducible polynomials in $\Fq[x]$ having even degrees.
Then $\Fq(x,\sqrt{f_1},\dots,\sqrt{f_6})$
is an unramified Galois extension of $\Fq(x,\sqrt{f})$, with Galois group
$(\Z/2\Z)^5$, in which both places of $\Fq(x,\sqrt{f})$ which contain $1/x$
are completely split.  It follows that~$A^{+}(q)>0$.

For odd $q$ it remains only to prove the second statement for $q \gg 0$.
Let $\alpha_1,\dots,\alpha_{r+1}$ be distinct elements of~$\Fq$,
and let $f(x)=(x-\alpha_1)\dots (x-\alpha_{r+1})$.
Suppose $r$ is even.  Then
$F:=\Fq(\sqrt{x-\alpha_1},\dots,\sqrt{x-\alpha_{r+1}})$
is an unramified Galois extension of $\Fq(x,\sqrt{f})$, with Galois
group~$(\Z/2\Z)^r$.
By Hurwitz, $\Fq(x,\sqrt{f})$ has genus $g:=r/2$ and $F$ has
genus $1+2^r(r/2-1)$.  By the Weil lower bound,
the number of rational places
$N$ on $F$ satisfies
$$N\geq q+1-2\sqrt{q}\cdot 2^r(r/2-1).$$
Exactly $N/2^r$ rational places of $\Fq(x,\sqrt{f})$ split
completely in~$F$.
Choose $r$ as the even integer nearest $(\log_2 q)/3$.
For $q \gg 0$, we have $r \ge 4$ and $N/2^r \ge (r-2)^2/4-1 \ge 0$,
so it is possible to choose a subset $S$ of these rational places
with $\#S=\lfloor (r-2)^2/4-1 \rfloor$.
By the lemma, the $(S,2)$-class field tower is infinite,
so a proof similar to that of Lemma~\ref{lemma-cft} shows that
    $$A^+(q) \ge \frac{\#S}{g-1} = \frac{r}{2} + O(1) = c' \log q + O(1)$$
as $q \rightarrow \infty$, for some $c'>0$.

The case of even $q$ can be treated in a completely analogous manner, using
Artin-Schreier covers.
%
%
\end{proof}

The proof of the theorem shows that there exists $c>0$ such that for
every~$q$, there exists a curve $C$ over $\Fq$ with $g(C)>1$ and a set
$S$ of rational places on $\Fq(C)$ such that $\#S>(c\log q)\cdot(g(C)-1)$
and such that the $(S,2)$-class field tower of $\Fq(C)$ is infinite.
As in Section~\ref{sec-cft}, it follows that for every $q$ we
have $A^{-}(q) \geq (c/4)\log q$.
We do not know whether this bound can be improved significantly
when $q$ is prime.



\end{document}